\documentclass[11pt]{amsart}
\usepackage[utf8]{inputenc}
\usepackage[T1]{fontenc}
\usepackage{amssymb,amsthm}
\usepackage{xcolor}
\usepackage[colorlinks=true,linkcolor=blue!50!black,citecolor=blue!50!black,urlcolor=blue!50!black]{hyperref}

\theoremstyle{plain}
\newtheorem{theorem}{Theorem}[section]
\newtheorem{lemma}[theorem]{Lemma}
\newtheorem{corollary}[theorem]{Corollary}
\newtheorem{proposition}[theorem]{Proposition}
\theoremstyle{definition}
\newtheorem{definition}[theorem]{Definition}
\newtheorem{example}[theorem]{Example}
\theoremstyle{remark}
\newtheorem{remark}[theorem]{Remark}

\newcommand{\fg}{\mathfrak g}

\newcommand{\fa}{\mathfrak a}
\newcommand{\fp}{\mathfrak p}
\newcommand{\R}{\mathbb R}
\newcommand{\Z}{\mathbb Z}
\newcommand{\ip}[2]{\langle #1,#2\rangle}
\newcommand{\conv}{\operatorname{conv}}
\newcommand{\Ad}{\operatorname{Ad}}
\newcommand{\spn}{\operatorname{span}}
\newcommand{\reg}{\mathrm{reg}}

\title[Eaton Triples, Polar Representations, and Lie Triple Systems]{Eaton Triples, Polar Representations, \\and Lie Triple Systems}
\author{Gabriel Larotonda}
\address{Instituto Argentino de Matem\'atica (CONICET) and Facultad de Ciencias Exactas y Naturales, Universidad de Buenos Aires, Argentina}
\email{glaroton@dm.uba.ar}
\thanks{Corresponding author: Gabriel Larotonda.}

\author{Tin-Yau Tam}
\address{Department of Mathematics and Statistics, University of Nevada, Reno, Nevada, USA}
\email{ttam@unr.edu}


\date{November 7, 2026}

\begin{document}
\begin{abstract}
An Eaton triple consists of a finite dimensional real inner product space, a closed subgroup of its orthogonal group, and a closed convex cone that meets every orbit and satisfies an extremal inner product inequality. In the literature on group majorization and invariant norms it is a standing hypothesis that the triple admits a reduced triple, obtained by restriction to the span of the cone. We prove that this is always the case: the axioms force the span of the cone to be a section of a polar representation, and Gichev's orbit intersection theorem applies. The reduced group is a finite reflection group, crystallographic when the group is connected; Eaton cones are its closed chambers, unique up to the group action; and a compact linear group admits an Eaton cone if and only if its representation is polar with a finite reflection group as Weyl group. For connected groups we construct a Lie triple system on the space, determined by the inner product and the orbits alone, for which the cone is a Weyl chamber; unlike Dadok's orbit equivalence with a symmetric space, it involves no choices.
\end{abstract}
\maketitle

\section{Introduction}

Let $V$ be a finite dimensional real inner product space and $G\subset \mathrm{O}(V)$ a closed subgroup of its orthogonal group, hence a compact Lie group, not necessarily connected. Following Eaton \cite[Ch.~6, assumption (A.I), p.~155]{Eaton}, a triple $(V,G,F)$ is an \emph{Eaton triple} if $F\subset V$ is a nonempty closed convex cone such that
\begin{itemize}
\item[(A1)] $Gx\cap F\neq\emptyset$ for every $x\in V$;
\item[(A2)] $\max_{g\in G}\ip{x}{gy}=\ip{x}{y}$ for all $x,y\in F$.
\end{itemize}
Conditions (A1) and (A2) imply that $Gx\cap F$ is a single point, denoted $\psi_F(x)$. The map $\psi_F:V\to F$ is the normal map of a normal decomposition system in the sense of Lewis \cite{Lewis}; in finite dimension the two notions are essentially equivalent, see \cite[Remark 1.5]{DolberryTam}. Examples include the symmetric group acting on $\R^n$, with $F$ the cone of decreasing vectors; the unitary group acting by conjugation on Hermitian matrices, with $F$ the decreasing diagonal matrices, where (A2) is Ky Fan's inequality; $\mathrm{U}(p)\times \mathrm{U}(q)$ acting on complex $p\times q$ matrices, with the inner product $\ip{x}{y}=\operatorname{Re}\operatorname{tr}(x^*y)$ and $F$ the real diagonal matrices with nonnegative decreasing entries, where (A2) is von Neumann's trace inequality; and the adjoint action of a compact Lie group on its Lie algebra, with $F$ a closed Weyl chamber, where (A2) follows from two results of Kostant (Section~\ref{sec:biinvariant}). The interplay between group majorization and reflection groups goes back to Eaton and Perlman \cite{EatonPerlman} and Steerneman \cite{Steerneman}.

Set $W=\spn F$, let
\[
N_G(W)=\{g\in G: gW=W\},\quad  Z_G(W)=\{g\in G: g|_W=\mathrm{id}\},
\]
 and
\[
H=\{g|_W:\ g\in N_G(W)\}\simeq N_G(W)/Z_G(W)\subset \mathrm{O}(W).
\]
Following Niezgoda \cite{Niezgoda} and Tam \cite{Tam}, $(W,H,F)$ is the \emph{reduced triple} of $(V,G,F)$ if it is itself an Eaton triple. It exists in all the examples above, with $H$ a crystallographic Weyl group, and is a standing hypothesis in the theory of $G$-invariant norms and group majorization on Eaton triples \cite{Tam,TamHill,NiemiecTam,DolberryTam}. Niezgoda \cite[Theorem 3.2]{Niezgoda} characterized its existence by a Schur type inequality and by the coincidence of the $G$- and $H$-majorizations on $W$ (Theorem~\ref{thm:niezgoda} below), but his theorem does not assert that it always exists. For finite $G$ the question was settled by Steerneman \cite[Theorem 4.1]{Steerneman}, see also \cite[p.~817]{Lewis2000}: $G$ is then a reflection group, $F$ is the closure of a chamber, and $(V,G,F)$ is its own reduced triple.

The main results are Theorems~\ref{thm:main}, \ref{thm:connected}, \ref{thm:rigidity}, \ref{thm:char}, and \ref{thm:canonical-triple}. Theorem~\ref{thm:main} proves the existence of the reduced triple for every Eaton triple. Theorem~\ref{thm:connected} describes the reduced group for connected $G$: it is a crystallographic Weyl group, as Palais and Terng observed \cite[Corollary 4.7]{PT}, and its fixed subspace is $V^G$. Theorem~\ref{thm:rigidity}, the rigidity of Eaton cones, rests on Theorem~\ref{thm:chamber}, the description of the cone in Steerneman's theorem applied to the finite group $H$. Theorem~\ref{thm:char} characterizes the representations admitting Eaton cones, and Theorem~\ref{thm:canonical-triple} constructs a canonical Lie triple system for connected groups.

The key observation behind Theorem~\ref{thm:main} is that the axioms force polarity: (A1) and the infinitesimal form of (A2) imply that $W$ is a section for $G$ (Lemma~\ref{lem:orth}). This places Eaton triples, which come from matrix analysis and group majorization, in the framework of polar representations. Gichev's orbit intersection theorem \cite[Proposition 2]{Gichev} then gives $Gx\cap W=Hx$, and Niezgoda's theorem gives the rest. Hence the existence of a reduced triple can be dropped from the hypotheses in \cite{Tam,TamHill,NiemiecTam,DolberryTam}, for instance from the description of $G$-invariant convex functions and norms by their restrictions to $W$ \cite[Theorems 4, 6 and 7]{TamHill} and from the subdifferential formula for a $G$-invariant convex function \cite[Theorem 4.9]{DolberryTam}. Theorem~\ref{thm:char} also rests on a theorem of Gichev. The conclusion of Theorem~\ref{thm:connected} fails for disconnected groups, finite or infinite (Section~\ref{sec:examples}).

A polar representation may be the isotropy representation of several symmetric spaces with the same orbits ($\mathrm{SO}(2m)$ and $\mathrm{U}(m)$ on $\R^{2m}$ is the simplest case), so Dadok's orbit equivalence involves a choice. We use the bracket construction of Cartan and Kostant, in the form given by Eschenburg and Heintze \cite[Sections~3--4, especially Lemma~4.1]{EH}, applied not to $G$ but to the largest connected group $\widehat G$ with the same orbits. The resulting \emph{trace-normalized} structure is canonical: a linear isometry between two Eaton triples with connected groups that maps orbits onto orbits is an isomorphism of their Lie triple systems (Theorem~\ref{thm:canonical-triple}). Other connected orbit-equivalent isotropy groups, with their trace metrics, give a different structure only on irreducible summands of cohomogeneity one (Remark~\ref{rem:noncanonical-lie-triple}).

Niezgoda \cite{Niezgoda2002} studied the structure of a class of Eaton triples from a different point of view; our emphasis is on the geometry of the group action, its polar structure and Weyl group, and the relation with symmetric spaces.

Earlier, Tam \cite[Theorems~2 and~7]{Tam1998} extended Ky Fan's eigenvalue and singular-value inequalities for sums of matrices to the isotropy representations of noncompact real semisimple Lie algebras and to the adjoint representations of compact connected Lie groups, using the partial order defined by the dual of a Weyl chamber. For elements of the closed chamber, this order is equivalent to majorization under the corresponding finite Weyl group. Corollary~\ref{cor:normal-majorization} recovers these inequalities in the Eaton-triple framework and establishes the same Weyl-group majorization inequality for arbitrary finite-dimensional Eaton triples, without assuming connectedness or the prior existence of a reduced triple.

Recent applications of Eaton triples to matrix inequalities include Niezgoda's refinement of von Neumann-type inequalities through product Eaton triples \cite{Niezgoda2023} and his interpolation of generalized Ky Fan inequalities by means of gradients of $G$-increasing functions \cite{Niezgoda2024}. The latter work uses the $G$-majorization inequality $\psi_F(x+y)\preceq_G\psi_F(x)+\psi_F(y)$. Our Corollary~\ref{cor:normal-majorization} shows that this inequality can be formulated in terms of the finite reflection group $H$, without a separate hypothesis that a reduced triple exists.

Eaton triples also induce Fan--Theobald--von Neumann (FTvN) systems. Gowda and Jeong \cite[Section~10]{GowdaJeong2023} studied reduced FTvN systems and established corresponding Lidskii-type and majorization results under a reduction hypothesis. For the FTvN system induced by an Eaton triple, Theorem~\ref{thm:main} ensures that this hypothesis is always satisfied, as explained after Corollary~\ref{cor:normal-majorization}. Further developments in the more general setting of semi-FTvN systems include the Lie-algebraic commutativity and optimization results of Gowda and Sossa \cite{GowdaSossa2026} and the work of Jeong and Gowda \cite{JeongGowda2026} on majorization, automorphism groups, and eigenvalue maps.

Section~\ref{sec:prelim} collects the facts on polar representations and reflection groups that we need. Theorem~\ref{thm:main} is proved in Section~\ref{sec:existence}, Theorems~\ref{thm:connected}, \ref{thm:chamber}, and \ref{thm:rigidity} in Section~\ref{sec:weyl}, Theorem~\ref{thm:char} in Section~\ref{sec:char} and Theorem~\ref{thm:canonical-triple} in Section~\ref{sec:lie-triple}. Section~\ref{sec:examples} contains the examples.

\section{Preliminaries}\label{sec:prelim}

Throughout, $(V,G,F)$ is an Eaton triple, and $W$, $N_G(W)$, $Z_G(W)$ and $H$ are as in the introduction. Let $G^0$ be the identity component of $G$ and $\fg\subset\mathfrak{so}(V)$ its Lie algebra. For $v\in V$, the tangent space to the orbit $Gv$ at $v$ is
\[
T_v(Gv)=\fg v:=\{\xi v:\xi\in\fg\}.
\]
Indeed, the differential at the identity of the orbit map $g\mapsto gv$ sends $\xi\in\fg$ to
\[
\left.\frac{d}{dt}\right|_{t=0}e^{t\xi}v=\xi v,
\]
and its image is $T_v(Gv)$. This formula also holds when $G$ is disconnected, since the orbits $Gv$ and $G^0v$ have the same tangent space at $v$. We write $\nu_v=(\fg v)^\perp$ for the normal space.

We recall Niezgoda's characterization of the reduced triple.

\begin{theorem}[Niezgoda {\cite[Theorem 3.2]{Niezgoda}}; see {\cite[Theorem 1.1]{DolberryTam}}]\label{thm:niezgoda}
Let $(V,G,F)$ be an Eaton triple and let $\pi:V\to W$ be the orthogonal projection. The following are equivalent:
\begin{enumerate}
\item for $x,y\in W$, $y\in\conv(Gx)$ if and only if $y\in\conv(Hx)$;
\item $(W,H,F)$ is an Eaton triple, i.e. the reduced triple of $(V,G,F)$ exists;
\item $\pi(x)\in\conv\big(H\psi_F(x)\big)$ for every $x\in V$.
\end{enumerate}
If they hold, $H$ is a finite reflection group.
\end{theorem}

Note that (A2) for $(W,H,F)$ is automatic: the elements of $H$ are restrictions of elements of $G$, so for $x,y\in F$ we have
\[
\max_{h\in H}\ip{x}{hy}\le\max_{g\in G}\ip{x}{gy}=\ip{x}{y},
\]
and $h=\mathrm{id}$ gives equality. Hence the existence of the reduced triple amounts to (A1) for $H$: for every $w\in W$ there is $n\in N_G(W)$ with $nw\in F$.

We use Gichev's definition of polar representation \cite{Gichev}, which does not assume connectedness; polar representations of connected groups were introduced by Dadok \cite{Dadok} and studied by Palais and Terng \cite{PT}. The orthogonal representation of the compact group $G$ on $V$ is \emph{polar} if there is a linear subspace $\fa\subset V$, a \emph{section} or \emph{Cartan subspace}, such that
\begin{itemize}
\item[(A)] every $G$-orbit meets $\fa$;
\item[(B)] $\fg u\perp\fa$ for every $u\in\fa$.
\end{itemize}
The \emph{Weyl group} of the section is
\[
\mathcal W(\fa)=\{g|_\fa: g\in G,\ g\fa=\fa\}\simeq N_G(\fa)/Z_G(\fa).
\]
The Weyl group $\mathcal W(\fa)$ is finite: if $\xi$ is in the Lie algebra of $N_G(\fa)$ then $\xi\fa\subset\fa$, and $\xi\fa\perp\fa$ by (B), so $\xi|_\fa=0$; hence the identity component of $N_G(\fa)$ acts trivially on $\fa$, and the compact Lie group $N_G(\fa)$ has finitely many components. By Theorem~\ref{thm:gichev-prop} below, $\fa=\nu_u$ for generic $u\in\fa$, and a section meets every orbit in finitely many points. A point is \emph{regular} if its isotropy group in $G$ is minimal; regular points have conjugate isotropy groups, hence orbits of the same dimension, and $V_\reg$ denotes the set of regular points.

\begin{theorem}[Gichev {\cite[Proposition 2]{Gichev}}]\label{thm:gichev-prop}
Let the representation of the compact group $G$ be polar with section $\fa$ and Weyl group $\mathcal W(\fa)$, and let $\pi:V\to\fa$ be the orthogonal projection. Then $\fa\cap V_\reg$ is open and dense in $\fa$, $\fa=\nu_v$ for $v\in\fa\cap V_\reg$, and for every $a\in\fa$
\[
Ga\cap\fa=\mathcal W(\fa)\,a\qquad\text{and}\qquad \pi\big(\conv(Ga)\big)=\conv\big(\mathcal W(\fa)\,a\big).
\]
\end{theorem}

\begin{theorem}[Gichev {\cite[Theorem 1]{Gichev}}]\label{thm:gichev-SP}
A compact linear group $G\subset \mathrm{O}(V)$ has the property that the family $\{\conv(Gv)\}_{v\in V}$ is a semigroup under Minkowski addition if and only if its representation is polar and its Weyl group is a finite reflection group. If $G$ is connected, this holds if and only if the representation is polar.
\end{theorem}

The last assertion is the theorem of Terng and of Palais and Terng that the Weyl group of a polar representation of a connected group is a finite reflection group (see the introduction of \cite{Gichev}). In fact it is a Weyl group in the Lie theoretic sense \cite[Corollary 4.7]{PT}; this follows from Dadok's orbit equivalence theorem below, see also the remark following \cite[Proposition~6]{Dadok}.

Two orthogonal representations are \emph{orbit equivalent} if a linear isometry between their spaces maps the orbits of one onto those of the other. The \emph{cohomogeneity} of a representation is the codimension of a regular orbit; for a polar representation it is the dimension of a section, by Theorem~\ref{thm:gichev-prop}. Following \cite[Introduction]{EH}, an \emph{$s$-representation} is the isotropy representation, on the tangent space $\fp$ at the base point, of a Riemannian symmetric space without Euclidean factor, restricted to the identity component of the isotropy group. The maximal abelian subspaces of $\fp$ are sections \cite[Section~4]{PT}, so its \emph{rank}, the dimension of a maximal abelian subspace of $\fp$, is its cohomogeneity. Dadok's orbit equivalence theorem \cite[Proposition~6]{Dadok} states that every polar representation of a connected compact Lie group is orbit equivalent to the isotropy representation of a Riemannian symmetric space; see \cite[Introduction and Main Theorem]{EH} for a conceptual approach in irreducible cohomogeneity at least three and a discussion of the lower-cohomogeneity cases.

A \emph{finite reflection group} is a finite subgroup of an orthogonal group generated by reflections (also called a finite Coxeter group in its reflection representation); we reserve \emph{crystallographic Weyl group} for the crystallographic case.

Finally, let $H\subset \mathrm{O}(W)$ be a finite reflection group, let $W_0=W^H$ be its fixed subspace and let $\mathcal A$ be the arrangement of its reflecting hyperplanes, all of which contain $W_0$; the chambers are the connected components of the complement of $\mathcal A$. We need the following facts, see \cite[Ch.~V, \S3, Th.~1, Th.~2 and Prop.~2]{Bourbaki46} and Humphreys \cite[Sections 1.12 and 2.8]{Humphreys}: $H$ acts simply transitively on the chambers; the closure $\overline C$ of a chamber is a fundamental domain and a closed convex cone with $\spn\overline C=W$ and $W_0\subset\overline C$; the isotropy group of a point is generated by the reflections in the hyperplanes containing it; and $\ip{x}{hy}\le\ip{x}{y}$ for $x,y\in\overline C$ and $h\in H$, because $y-hy$ is a nonnegative combination of the roots that are positive on $C$. The group $H$ is \emph{crystallographic} if it is the Weyl group of a crystallographic root system in $W_0^\perp$; for a dihedral group $I_2(m)$ this happens exactly when $m\in\{2,3,4,6\}$.

\section{Existence of the reduced triple}\label{sec:existence}

\begin{lemma}[Orthogonality]\label{lem:orth}
We have $\ip{x}{\xi y}=0$ for every $\xi\in\fg$ and $x,y\in W$. In particular $W\subset\nu_w$ for every $w\in W$, and $gW\subset\nu_{gw}$ for every $g\in G$ and $w\in W$.
\end{lemma}
\begin{proof}
Let $x,y\in F$ and $\xi\in\fg$. By (A2) the differentiable function $t\mapsto\ip{x}{e^{t\xi}y}$ attains its maximum at $t=0$, so its derivative there, $\ip{x}{\xi y}$, vanishes. Since $F$ is a convex cone, $W=\spn F=F-F$ (collect the terms with positive and with negative coefficients), and the identity extends to $W$ by bilinearity.

For the last assertion, let $g\in G$, $w\in W$ and $\xi\in\fg$. Since $G$ is a matrix group, $g^{-1}\xi g\in\fg$, and
\[
\xi(gw)=g\bigl((g^{-1}\xi g)w\bigr).
\]
By the first assertion $(g^{-1}\xi g)w\perp W$, and $g$ is orthogonal, so $\xi(gw)\perp gW$. Since $\xi\in\fg$ is arbitrary, $gW\subset\nu_{gw}$.
\end{proof}

Together with (A1), the lemma says that $W$ is a section in the sense of Section~\ref{sec:prelim}: the representation of $G$ is polar. Moreover $V^G\subset F\subset W$, where $V^G$ is the space of $G$-fixed vectors: a fixed vector is its own orbit and meets $F$ by (A1).

\begin{theorem}\label{thm:main}
Every Eaton triple $(V,G,F)$ has a reduced triple $(W,H,F)$, and $H$ is finite. Moreover $Gw\cap W=N_G(W)\,w=Hw$ for every $w\in W$.
\end{theorem}
\begin{proof}
Since $W$ is a section, Theorem~\ref{thm:gichev-prop} gives $Gw\cap W=\mathcal W(W)\,w=N_G(W)\,w$ for $w\in W$, and $\mathcal W(W)=H$ by definition. Given $w\in W$, (A1) gives $g\in G$ with $gw\in F\subset W$, hence $gw\in Hw$; this is (A1) for $H$, and (A2) was checked in Section~\ref{sec:prelim}.

Finally, if $\xi$ lies in the Lie algebra of $N_G(W)$ then $\xi W\subset W$, and $\xi W\perp W$ by Lemma~\ref{lem:orth}, so $\xi|_W=0$. Hence the identity component of $N_G(W)$ acts trivially on $W$, and $H$ is finite because the compact Lie group $N_G(W)$ has finitely many connected components.
\end{proof}

\begin{corollary}\label{cor:niezgoda}
For every Eaton triple, $H$ is a finite reflection group, $\pi(x)\in\conv(H\psi_F(x))$ for every $x\in V$, and for $x,y\in W$ one has $y\in\conv(Gx)$ if and only if $y\in\conv(Hx)$.
\end{corollary}
\begin{proof}
Theorem~\ref{thm:main} is condition (2) of Theorem~\ref{thm:niezgoda}, and the rest follows from the equivalences there.
\end{proof}

\begin{corollary}[Reduction of group majorization]\label{cor:normal-majorization}
For an Eaton triple $(V,G,F)$, write $y\preceq_G x$ if
$y\in\conv(Gx)$, and define $\preceq_H$ analogously on $W$.
Then for all $x,y\in V$,
\[
y\preceq_G x\quad\Longleftrightarrow\quad
\psi_F(y)\preceq_H\psi_F(x).
\]
Moreover, one recovers the Fan-type inequality
\cite[Theorem~10]{Tam2000} without a separate reduced-triple hypothesis:
\[
\psi_F(x+y)\preceq_H\psi_F(x)+\psi_F(y)
\qquad (x,y\in V).
\]
\end{corollary}
\begin{proof}
Every vector lies in the $G$-orbit of its image under $\psi_F$, so $y\preceq_G x$ if and only if $\psi_F(y)\preceq_G\psi_F(x)$, which is equivalent to $\psi_F(y)\preceq_H\psi_F(x)$ by Corollary~\ref{cor:niezgoda}. For the last assertion, (A2) gives, for every $z\in F$,
\begin{align*}
\ip{z}{\psi_F(x+y)}
 &=\max_{g\in G}\ip{z}{g(x+y)}\\
 &\le \max_{g\in G}\ip{z}{gx}+\max_{g\in G}\ip{z}{gy}\\
 &=\ip{z}{\psi_F(x)}+\ip{z}{\psi_F(y)}.
\end{align*}
To see that these inequalities imply $H$-majorization, let $a,b\in F$ satisfy $\ip{z}{a}\le\ip{z}{b}$ for all $z\in F$, let $u\in W$, and choose $h\in H$ with $hu\in F$ by (A1). Then (A2) gives
\[
\ip{u}{a}=\ip{hu}{ha}\le\ip{hu}{a}
\le\ip{hu}{b}=\max_{k\in H}\ip{u}{kb}.
\]
Here the first inequality is (A2) for the pair $hu,a\in F$, and the second is our assumption with $z=hu$. By the support-function characterization of convex hulls, $a\in\conv(Hb)$, that is, $a\preceq_H b$. With $a=\psi_F(x+y)$ and $b=\psi_F(x)+\psi_F(y)$ this proves the claim.
\end{proof}

For the isotropy representations of noncompact symmetric spaces and the adjoint representations of compact connected Lie groups, the last inequality specializes to the Weyl-chamber inequalities established in \cite[Theorems~2 and~7]{Tam1998}, respectively. It is also the finite-reflection-group formulation of the generalized Ky Fan inequality used in \cite{Niezgoda2024}. Moreover, let $\lambda=\psi_F:V\to W$ and $\mu=\psi_F|_W:W\to W$. By Theorem~\ref{thm:main}, $(W,H,F)$ is an Eaton triple, so $(W,W,\mu)$ is a reduced FTvN system of $(V,W,\lambda)$ in the sense of \cite[Definition~10.1]{GowdaJeong2023}: indeed, $\mu\circ\lambda=\lambda$ and $\operatorname{ran}\mu=\operatorname{ran}\lambda=F$. Hence \cite[Theorem~10.3 and Corollary~10.4]{GowdaJeong2023} apply automatically to FTvN systems induced by finite-dimensional Eaton triples.

Connectedness is not used, but the components of $G$ cannot be ignored: Example~\ref{ex:swap} gives an Eaton triple $(V,G,F)$ for which $(V,G^0,F)$ is not one.

\section{Weyl groups and the rigidity of Eaton cones}\label{sec:weyl}

By Corollary~\ref{cor:niezgoda} the group $H$ is a finite reflection group. For connected $G$ it is a Weyl group in the Lie theoretic sense: this is the first assertion of the next theorem, and it is \cite[Corollary 4.7]{PT}, a consequence of Dadok's orbit equivalence theorem. We include a proof because it identifies $H$ with the restricted Weyl group of the symmetric space provided by Dadok's theorem, which gives the second assertion.

\begin{theorem}\label{thm:connected}
Let $(V,G,F)$ be an Eaton triple with $G$ connected. Then $H$ is a crystallographic Weyl group, and its fixed subspace $W^H$ is $V^G$.
\end{theorem}
\begin{proof}
By Section~\ref{sec:existence}, $W$ is a section for the polar representation of $G$. Let $V'=(V^G)^\perp$. Since $V^G\subset W$ and $G$ fixes $V^G$ pointwise,
\[
V=V^G\oplus V',\qquad W=V^G\oplus W',\qquad W'=W\cap V'.
\]
The subspace $W'$ is a section for the action on $V'$: an orbit contained in the $G$-invariant subspace $V'$ meets $W$ by (A1), hence meets $W\cap V'=W'$, and condition (B) is inherited from $W$.

By Dadok's orbit equivalence theorem, we may identify $V'$ isometrically with the space $\fp$ of the isotropy representation of a Riemannian symmetric space, with connected isotropy group $K$, so that the $G$-orbits and the $K$-orbits coincide. Since $G$ has no nonzero fixed vectors in $V'$, the symmetric space has no Euclidean factor, and we may take it of noncompact type. As the orbits are the same, $W'$ is also a section for $K$, and we claim that it is a maximal abelian subspace of $\fp$. Indeed, by Theorem~\ref{thm:gichev-prop} applied to the action of $K$ on $V'$, there is a $K$-regular point $u\in W'$, and $W'$ is the normal space to $Ku$ in $V'$. The point $u$ lies in a maximal abelian subspace $\fa$ of $\fp$, which is a section for $K$ (Section~\ref{sec:prelim}), so by the same theorem $\fa$ is this normal space, and $W'=\fa$. Extend the action of $K$ trivially to $V^G$. Then $G$ and $K$ have the same orbits on $V$, and $W$ is a section for both. Write
\[
H_G=N_G(W)/Z_G(W),\qquad H_K=N_K(W)/Z_K(W).
\]
By Theorem~\ref{thm:main}, $Gx\cap W=H_Gx$, and Theorem~\ref{thm:gichev-prop} applied to the section $W$ of $K$ gives $Kx\cap W=H_Kx$. Thus
\[
H_Gx=Gx\cap W=Kx\cap W=H_Kx\qquad(x\in W).
\]

We claim that $H_G=H_K$ as subgroups of $\mathrm{O}(W)$. Fix $h\in H_G$. For every $x\in W$ there is $k\in H_K$ with $hx=kx$, so
\[
W=\bigcup_{k\in H_K}\ker(h-k).
\]
Since $H_K$ is finite and a real vector space is not the union of finitely many proper linear subspaces, $\ker(h-k)=W$ for some $k\in H_K$. Hence $h=k$, and $H_G\subset H_K$. The reverse inclusion follows by interchanging the two groups.

The group $H_K$ is the Weyl group of the restricted root system, hence crystallographic, even when this system is nonreduced \cite[Ch.~VII]{Helgason}. Moreover, the restricted roots span the dual of $W'$: otherwise their common kernel would give a nonzero Euclidean factor of the symmetric space, contrary to the choice of noncompact type. Hence the vectors of $W$ fixed by $H_K$ are exactly those of $V^G$, and since $H=H_G=H_K$ we obtain $W^H=V^G$.
\end{proof}

Section~\ref{sec:examples} shows that the conclusion fails for disconnected $G$, finite or infinite, with $W=V$ or with $W\subsetneq V$.

\begin{theorem}[Eaton cones are closed chambers]\label{thm:chamber}
Let $H\subset \mathrm{O}(W)$ be a finite reflection group and let $F\subset W$ be a closed convex cone with $\spn F=W$ such that $(W,H,F)$ is an Eaton triple. Then $F$ is the closure of a chamber of $H$.
\end{theorem}

This is contained in Steerneman's theorem \cite[Theorem 4.1]{Steerneman}, which shows moreover that a finite group admitting an Eaton cone is a reflection group; we include the short proof for completeness.

\begin{proof}
If $W=\{0\}$, the assertion is immediate. Otherwise, since $F$ spans $W$, it contains $d=\dim W$ linearly independent vectors. Their nonnegative span is a full-dimensional simplicial cone contained in $F$, so $\operatorname{int}_W F\neq\varnothing$.

Let $s_\alpha\in H$ be the reflection in the hyperplane $\alpha^\perp$. For $x,y\in F$, (A2) gives
\[
0\le\ip{x}{y}-\ip{x}{s_\alpha y}=\frac{2\ip{x}{\alpha}\ip{y}{\alpha}}{\ip{\alpha}{\alpha}}.
\]
Hence $\ip{\cdot}{\alpha}$ does not change sign on $F$, and $F$ lies in one of the two closed half-spaces bounded by $\alpha^\perp$. Since $\spn F=W$, the functional $\ip{\cdot}{\alpha}$ cannot vanish identically on $F$. It therefore has a fixed, strictly nonzero sign on $\operatorname{int}_W F$: otherwise an open ball contained in $F$ and centered at a point of $\alpha^\perp$ would meet both sides of that hyperplane. Consequently $\operatorname{int}_W F$ is connected and disjoint from every reflecting hyperplane, so it lies in a chamber $C$. Since a closed convex set with nonempty interior is the closure of its interior, $F=\overline{\operatorname{int}_W F}\subset\overline C$. Conversely, let $c\in C$. By (A1) there is $h\in H$ with $hc\in F\subset\overline C$. Since $c$ lies on no reflecting hyperplane, neither does $hc$, so $hc\in C$ and $hC=C$. As $H$ acts simply transitively on the chambers, $h=\mathrm{id}$ and $c\in F$. Hence $F=\overline C$.
\end{proof}

\begin{theorem}[Rigidity]\label{thm:rigidity}
Let $G\subset \mathrm{O}(V)$ be closed and let $F_1,F_2\subset V$ be closed convex cones such that $(V,G,F_1)$ and $(V,G,F_2)$ are Eaton triples. Then there is $g\in G$ with $gF_1=F_2$. Consequently $\psi_{F_2}(x)=g\psi_{F_1}(x)$ for every $x\in V$, the reduced triples of the two Eaton triples are conjugate under $g$, and so are their groups $H_1$ and $H_2$.
\end{theorem}
\begin{proof}
Put $W_i=\spn F_i$. By Section~\ref{sec:existence}, $W_1$ and $W_2$ are sections for $G$; we first show that they are conjugate. Choose a regular point $w\in W_1$, which exists by Theorem~\ref{thm:gichev-prop}, and by (A1) for $F_2$ choose $g_0\in G$ with $g_0w\in F_2\subset W_2$. The point $g_0w$ is also regular, and Lemma~\ref{lem:orth}, applied to $(V,G,F_2)$, gives
\[
W_2\subset\nu_{g_0w}=g_0\nu_w=g_0W_1.
\]
Here $\nu_w=W_1$ by Theorem~\ref{thm:gichev-prop}. The same theorem gives $W_2=\nu_u$ for a regular point $u\in W_2$, and all regular orbits have the same dimension, so $\dim W_1=\dim W_2$ and $g_0W_1=W_2$.

Since $g_0\in G$, $(V,G,g_0F_1)$ is an Eaton triple, and the cones $g_0F_1$ and $F_2$ both span $W_2$. By Theorem~\ref{thm:main}, Corollary~\ref{cor:niezgoda} and Theorem~\ref{thm:chamber}, both are closed chambers of the finite reflection group $H=N_G(W_2)/Z_G(W_2)$, so $h(g_0F_1)=F_2$ for some $h\in H$. If $n\in N_G(W_2)$ restricts to $h$ on $W_2$, then $g=ng_0$ satisfies $gF_1=F_2$.

For $x\in V$, the point $g\psi_{F_1}(x)$ lies in $F_2\cap Gx=\{\psi_{F_2}(x)\}$, which proves the identity for the normal maps. Finally, conjugation by $g$ takes $N_G(W_1)$ and $Z_G(W_1)$ onto $N_G(W_2)$ and $Z_G(W_2)$ respectively, and therefore identifies the reduced triples and their groups.
\end{proof}

Theorems~\ref{thm:chamber} and \ref{thm:rigidity} give the claimed rigidity. In Theorem~\ref{thm:chamber} the triple $(W,H,F)$, with $W=\spn F$, is its own reduced triple, so $H$ is generated by reflections by Corollary~\ref{cor:niezgoda}.

\section{Examples}\label{sec:examples}

\begin{example}[The components matter]\label{ex:swap}
Let $G=\mathrm{SO}(2)\times\{\pm1\}$ act on $V=\R^2\oplus\R$ by
\[
(R,\varepsilon)(u,t)=(Ru,\varepsilon t).
\]
Set
\[
F=\{(r,0,s):\ r\ge0,\ s\ge0\}.
\]
Then $(V,G,F)$ is an Eaton triple: the orbit of $(u,t)\in V$ contains $(\|u\|,0,|t|)\in F$, which is (A1), and for $x=(r,0,s)$ and $y=(r',0,s')$ in $F$,
\[
\ip{x}{(R,\varepsilon)y}=rr'\ip{e_1}{Re_1}+\varepsilon ss'
\le rr'+ss'=\ip{x}{y},
\]
which is (A2). On the other hand, $G^0=\mathrm{SO}(2)\times\{1\}$ fixes the last coordinate, so the $G^0$-orbit of $(0,0,-1)$ does not meet $F$ and $(V,G^0,F)$ is not an Eaton triple.

Here $W=\spn\{e_1,e_3\}$, and $H\simeq\Z_2\times\Z_2$ acts on $W$ by independent sign changes; it is crystallographic, of type $A_1\times A_1$, so the conclusion of Theorem~\ref{thm:connected} holds for this disconnected group. The next two examples show that it can fail.
\end{example}

\begin{example}[A finite non-crystallographic example]\label{ex:I25}
Let $G=I_2(5)$ be the dihedral group of order $10$ generated by the reflections of $V=\R^2$ in two lines at angle $\pi/5$, and let $F$ be the closed sector between them. The ten sectors $gF$ tile the plane, which is (A1), and (A2) is the chamber inequality recalled in Section~\ref{sec:prelim}. Since $F$ has nonempty interior, $W=V$, $H=G$ and $(V,G,F)$ is its own reduced triple; by the criterion in Section~\ref{sec:prelim}, $H$ is a finite reflection group but not crystallographic.
\end{example}

\begin{example}[An infinite disconnected example with $W\subsetneq V$]\label{ex:SO2I25}
Let $D=I_2(5)$ act on $\R^2$ as above, with closed chamber $C$, and let
\[
V=\R^2\oplus\R^2,\qquad G=\mathrm{SO}(2)\times D,\qquad F=\R_{\ge0}e_1\oplus C .
\]
Then $G$ is compact and infinite, with ten components. For $(u,v)\in V$ choose $r\in \mathrm{SO}(2)$ with $ru=\|u\|e_1$ and $d\in D$ with $dv\in C$; then $(r,d)(u,v)\in F$, which is (A1). For $x=(ae_1,u)$ and $y=(be_1,v)$ in $F$ and $(r,d)\in G$,
\[
\ip{x}{(r,d)y}=ab\ip{e_1}{re_1}+\ip{u}{dv}\le ab+\ip{u}{v}=\ip{x}{y},
\]
by the chamber inequality for $D$, which is (A2). Here $W=\R e_1\oplus\R^2$ has dimension $3$, $N_G(W)=\{\pm I\}\times D$, $Z_G(W)$ is trivial, and
\[
H\simeq\Z_2\times I_2(5),
\]
a non-crystallographic finite reflection group. The identity $Gx\cap W=Hx$ is immediate: for $x=(ae_1,u)\in W$, $Gx\cap W=\{\pm ae_1\}\times Du=Hx$.
\end{example}

\begin{remark}[Dihedral groups]
The same construction with the dihedral group $I_2(m)$ of order $2m$ gives Eaton triples whose group $H$ is $I_2(m)$ or $\Z_2\times I_2(m)$, which are crystallographic exactly when $m\in\{2,3,4,6\}$. Thus $I_2(5)$ is the smallest dihedral group for which the conclusion of Theorem~\ref{thm:connected} fails.
\end{remark}

\begin{example}[A polar action with no Eaton cone]\label{ex:C3}
Let $G=C_3\subset \mathrm{SO}(2)$ be the group generated by the rotation $\rho$ of $\R^2$ by $2\pi/3$. Being finite, $G$ is polar with section $\R^2$ and Weyl group $C_3$, which is not generated by reflections. This action admits no Eaton cone. To see the angular obstruction to (A1), note that the $G$-orbit of a unit vector with direction $\theta$ consists of the three directions $\theta$, $\theta+2\pi/3$, and $\theta+4\pi/3$ (modulo $2\pi$). A closed convex cone that is the origin, a ray, or a line misses some of these orbits. If $F$ is a sector of aperture $\beta<2\pi/3$, its three rotations by $G$ cover at most $3\beta<2\pi$ of the unit circle, so (A1) again fails. Thus any cone satisfying (A1) is either the entire plane or a closed sector of aperture at least $2\pi/3$ (including a half-plane). In either case one can choose a nonzero $x\in F$ with $\rho x\in F$: take $x$ on the initial boundary ray of the sector, or arbitrary if $F=\R^2$. But then
\[
\max_{g\in G}\ip{x}{g\rho x}=\ip{x}{\rho^3x}=\|x\|^2>-\tfrac12\|x\|^2=\ip{x}{\rho x},
\]
which violates (A2). This is the case excluded by Theorem~\ref{thm:char} below.
\end{example}

\section{Which actions admit an Eaton cone}\label{sec:char}

\begin{theorem}\label{thm:char}
Let $G\subset \mathrm{O}(V)$ be closed. The following are equivalent.
\begin{enumerate}
\item There is a closed convex cone $F$ such that $(V,G,F)$ is an Eaton triple.
\item The representation of $G$ is polar and the Weyl group of some (equivalently, of every) section is a finite reflection group.
\item The family $\{\conv(Gv)\}_{v\in V}$ is a semigroup under Minkowski addition.
\end{enumerate}
When they hold, the Eaton cones are the closed chambers of the Weyl groups of the sections. If $G$ is connected, the three conditions are equivalent to the representation being polar.
\end{theorem}
\begin{proof}
(1)$\Rightarrow$(2). By Lemma~\ref{lem:orth} and (A1) the representation is polar with section $W$, and by Corollary~\ref{cor:niezgoda} its Weyl group $H$ is a finite reflection group. If $\fa_1,\fa_2$ are sections, $w\in\fa_1$ is regular and $gw\in\fa_2$, then $\fa_2\subset\nu_{gw}=g\fa_1$ by condition (B) and Theorem~\ref{thm:gichev-prop}, with equality by dimension, as in the proof of Theorem~\ref{thm:rigidity}. Hence every section is conjugate to $W$ under $G$, and its Weyl group to $H$.

(2)$\Rightarrow$(1). Let $\fa$ be a section with $\mathcal W(\fa)$ a finite reflection group, let $C$ be a chamber of $\mathcal W(\fa)$ and put $F=\overline C$. By (A) every $G$-orbit meets $\fa$, and every $\mathcal W(\fa)$-orbit in $\fa$ meets $\overline C$ since the chambers tile $\fa$; this is (A1). Let $\pi$ be the orthogonal projection onto $\fa$. For $x,y\in F$ and $g\in G$, since $x\in\fa$,
\begin{align*}
\ip{x}{gy}=\ip{x}{\pi(gy)}&\le\max\{\ip{x}{z}: z\in\pi(\conv(Gy))\}\\
&=\max_{w\in \mathcal W(\fa)}\ip{x}{wy}=\ip{x}{y},
\end{align*}
where the second equality is Theorem~\ref{thm:gichev-prop} and the last one the chamber inequality; this is (A2).

(2)$\Leftrightarrow$(3) is Theorem~\ref{thm:gichev-SP}, which also gives the last assertion. Finally, if $(V,G,F)$ is an Eaton triple, then $W=\spn F$ is a section by Lemma~\ref{lem:orth} and (A1), $(W,H,F)$ is an Eaton triple with $H=\mathcal W(W)$ a finite reflection group by Theorem~\ref{thm:main} and Corollary~\ref{cor:niezgoda}, and $F$ is a closed chamber of $H$ by Theorem~\ref{thm:chamber}. Conversely, by the proof of (2)$\Rightarrow$(1) every closed chamber of the Weyl group of any section is an Eaton cone, and all these cones are conjugate under $G$ by Theorem~\ref{thm:rigidity}.
\end{proof}

For connected $G$ the theorem says that the Eaton triples are exactly the polar representations with a chosen Weyl chamber; by Dadok's theorem these are, up to orbit equivalence, the isotropy representations of symmetric spaces with a closed Weyl chamber of the restricted root system, the Eaton triples considered by Lewis \cite[Propositions 3.3 and 3.4]{Lewis2000} and by Tam \cite{Tam}. For disconnected $G$ the class is larger, as the dihedral examples show, and the condition in (2) on the Weyl group is a genuine restriction (Example~\ref{ex:C3}).

\section{The canonical Lie triple system of a connected Eaton triple}\label{sec:lie-triple}

In Subsections~\ref{sec:orbit-hull}--\ref{sec:noncanonical} $G$ is connected. We show that $V$ carries a Lie triple product, given by an explicit formula and determined by the inner product and the orbits, for which $W$ is a maximal abelian subspace and $F$ a closed Weyl chamber. The construction is that of Eschenburg and Heintze (trace metric: \cite[Section~3]{EH}; bracket and Jacobi-identity criterion: \cite[Section~4, especially Lemma~4.1]{EH}), applied not to $G$ but to the largest connected group with the same orbits.

\subsection{The orbit hull}\label{sec:orbit-hull}

\begin{definition}
Let $\fg\subset\mathfrak{so}(V)$ be the Lie algebra of $G$. The \emph{orbit hull} of $\fg$ is
\[
\widehat\fg=\{A\in\mathfrak{so}(V):\ Ax\in\fg x\ \text{ for every }x\in V\},
\]
the space of linear vector fields on $V$ which are tangent to all the $G$-orbits. Let $\mathrm{O}(V)_G=\{g\in \mathrm{O}(V):\ g(Gx)=Gx\ \text{for every }x\in V\}$ and let $\widehat G$ be its identity component.
\end{definition}

\begin{remark}[What the orbit hull is]\label{rem:hull-meaning}
An element $A\in\widehat\fg$ satisfies $Ax=\xi x$ for every $x\in V$, with $\xi\in\fg$ depending on $x$. For instance, $\mathrm{U}(m)$ and $\mathrm{SO}(2m)$ have the same orbits on $\mathbb C^m=\R^{2m}$, the spheres, so $\widehat\fg=\mathfrak{so}(2m)$: the orbit hull does not retain the complex structure. In the same way, $G_2\subset \mathrm{SO}(7)$ gives $\widehat\fg=\mathfrak{so}(7)$. On the other hand, for the adjoint representation of a compact semisimple Lie algebra $\widehat\fg=\fg$ (Remark~\ref{rem:group-case}).
\end{remark}

\begin{lemma}\label{lem:hull}
$\mathrm{O}(V)_G$ is a closed subgroup of $\mathrm{O}(V)$ with Lie algebra $\widehat\fg$, and $\widehat G$ is the largest connected subgroup of $\mathrm{O}(V)$ with the same orbits as $G$. Moreover $(V,\widehat G,F)$ is an Eaton triple, and its group $N_{\widehat G}(W)/Z_{\widehat G}(W)$ coincides with $H$ as a subgroup of $\mathrm{O}(W)$. Finally, let 
\[
V=V^G\oplus V_1\oplus\dots\oplus V_k
\]
 be an orthogonal decomposition into $G$-invariant subspaces, with $V_1,\dots,V_k$ irreducible and nontrivial, and let $\widehat\fg_i$ and $\widehat G_i$ be the orbit hull and the corresponding group for the action of $G$ on $V_i$. Then every element of $\widehat\fg$ preserves each $V_i$ and vanishes on $V^G$, $\widehat\fg=\bigoplus_i\widehat\fg_i$, $\widehat G=\widehat G_1\times\dots\times\widehat G_k$ (acting trivially on $V^G$), and each $\widehat G_i$ acts on $V_i$ as an irreducible $s$-representation.
\end{lemma}
\begin{proof}
Each orbit $Gx$ is compact, hence closed, so its setwise stabilizer in $\mathrm{O}(V)$ is closed, and $\mathrm{O}(V)_G$, the intersection of these stabilizers, is a closed subgroup of $\mathrm{O}(V)$.

If $A\in\mathfrak{so}(V)$ and $e^{tA}\in\mathrm{O}(V)_G$ for every $t$, the curve $t\mapsto e^{tA}x$ lies in $Gx$, and differentiating at $t=0$ gives
\[
Ax\in T_x(Gx)=\fg x\qquad(x\in V).
\]
Conversely, if $Ax\in\fg x$ for every $x\in V$, the linear vector field $x\mapsto Ax$ is tangent to every $G$-orbit, and its integral curve through $x$ is $t\mapsto e^{tA}x$. Since the integral curves of a vector field tangent to a compact embedded submanifold stay in it, $e^{tA}x\in Gx$ for all $t$, and the same applied to $e^{-tA}$ gives $e^{tA}(Gx)=Gx$. Thus $e^{tA}\in\mathrm{O}(V)_G$, and the Lie algebra of $\mathrm{O}(V)_G$ is $\widehat\fg$.

Since $G$ is connected and preserves each $G$-orbit, $G\subset\widehat G$, and since $\widehat G$ preserves each $G$-orbit, $\widehat Gx=Gx$ for all $x\in V$. A connected subgroup with these orbits lies in $\mathrm{O}(V)_G$, hence in its identity component $\widehat G$; this proves maximality. Since the orbits agree, (A1) holds for $(V,\widehat G,F)$, and for $x,y\in F$
\[
\max_{g\in\widehat G}\ip{x}{gy}
=\max_{z\in\widehat Gy}\ip{x}{z}
=\max_{z\in Gy}\ip{x}{z}
=\ip{x}{y}.
\]
Thus (A2) holds. Theorem~\ref{thm:main}, applied to both triples, gives
\[
N_{\widehat G}(W)x=\widehat Gx\cap W=Gx\cap W=N_G(W)x\qquad(x\in W).
\]
The images of the two normalizers in $\mathrm{O}(W)$ are finite subgroups, the first containing the second since $N_G(W)\subset N_{\widehat G}(W)$, and their point orbits agree; by the finite-group argument in the proof of Theorem~\ref{thm:connected}, they coincide.

The representation of $G$ is polar by Section~\ref{sec:existence}. By Dadok's decomposition theorem \cite[Theorem~4]{Dadok}, for an orthogonal $G$-invariant splitting $V=E_1\oplus E_2$ the section splits accordingly and the orbit of a point of the section is the product of its orbits in $E_1$ and $E_2$; moving an arbitrary $x=x_1+x_2$ into a section, $Gx=Gx_1\times Gx_2$. This fails for reducible actions in general. Iterating over the invariant summands and noting that $G$ fixes $V^G$ pointwise gives, for $x=x_0+\sum_i x_i$ with $x_0\in V^G$ and $x_i\in V_i$,
\[
Gx=\Big\{x_0+\sum_i y_i:\ y_i\in Gx_i\Big\}.
\]
The components $x_i$ thus move independently. Taking tangent spaces at $x$, and regarding each $T_{x_i}(Gx_i)$ as a subspace of $V_i$, we get
\[
\fg x=T_x(Gx)=\bigoplus_i T_{x_i}(Gx_i)=\bigoplus_i\fg x_i.
\]
The actions on the $V_i$ are polar by the same theorem; see also \cite[Proposition 3.4]{GT} and \cite{Bergmann} for related results.

If $A\in\widehat\fg$ and $x\in V_i$ then $Ax\in\fg x\subset V_i$, so $A$ preserves $V_i$; if $x\in V^G$ then $Ax\in\fg x=0$. The restriction of $A$ to $V_i$ lies in $\widehat\fg_i$, whence $\widehat\fg\subset\bigoplus_i\widehat\fg_i$. Conversely, if $A_i\in\widehat\fg_i$ for each $i$, the operator $A=\bigoplus_i A_i$, acting as zero on $V^G$, satisfies
\[
Ax=\sum_i A_i x_i\in\bigoplus_i\fg x_i=\fg x.
\]
Thus $A\in\widehat\fg$, and $\widehat\fg=\bigoplus_i\widehat\fg_i$. The group $\widehat G_1\times\dots\times\widehat G_k$, acting trivially on $V^G$, is a connected subgroup of $\mathrm{O}(V)$ with Lie algebra $\bigoplus_i\widehat\fg_i=\widehat\fg$, hence equals $\widehat G$.

Finally, by the first part of the proof $\widehat G_i$ is the largest connected group with the orbits of $G$ on $V_i$. It acts irreducibly, since it contains the image of $G$ on the irreducible space $V_i$, and its action is polar, since it has the orbits of $G$. In cohomogeneity one the orbits are spheres, so $\widehat G_i=\mathrm{SO}(V_i)$, the isotropy group of real hyperbolic space. Suppose that the cohomogeneity is at least two. By Dadok's theorem \cite[Proposition~6]{Dadok}, after an isometric identification of $V_i$ with the tangent space of a symmetric space there is an $s$-representation $K_i\subset\mathrm{O}(V_i)$ with the orbits of $G$, $K_i$ being the connected isotropy group. Groups with the same orbits have the same invariant subspaces (an invariant subspace is a union of orbits), so $K_i$ is irreducible, and its rank is the cohomogeneity of $V_i$. Since $K_i$ is connected and has the orbits of $G$, we have $K_i\subset\widehat G_i$. By \cite[Proposition~5.1]{EH}, an irreducible $s$-representation of rank at least two is maximal among connected subgroups of $\mathrm{O}(V_i)$ having the same orbits. Thus $K_i=\widehat G_i$.
\end{proof}

\subsection{The canonical Lie triple system}\label{sec:canonical}

Following Cartan and Kostant, we extend the Lie bracket of $\widehat\fg$ to $L=\widehat\fg\oplus V$ so that $(L,\widehat\fg)$ becomes a symmetric Lie algebra: the bracket of $\widehat\fg$ with $V$ is its action on $V$, and the bracket of two vectors of $V$ is defined through an invariant inner product on $\widehat\fg$. The Jacobi identity is what has to be verified; for this we use the trace metric and Lemma~4.1 of \cite{EH}. Since the orbit hull $\widehat\fg$ replaces $\fg$, the construction depends only on the orbits.

The \emph{trace metric} on $\widehat\fg$ is
\[
\tau(\xi,\eta)=-\operatorname{tr}(\operatorname{ad}\xi\circ\operatorname{ad}\eta)-\operatorname{tr}(\xi\eta)\qquad(\xi,\eta\in\widehat\fg),
\]
the first trace taken on $\widehat\fg$ and the second on $V$. It is $\Ad(\widehat G)$-invariant and positive definite: both terms are nonnegative and the second is the Hilbert--Schmidt norm on $\mathfrak{so}(V)$. The first term is explained by the symmetric case: when $\widehat G$ acts on $V$ as the isotropy representation of a symmetric space, $-\tau$ is the restriction to the isotropy algebra of the Killing form of the full algebra. On an irreducible summand $V_i$ the Killing form is a negative multiple of the given inner product, so the bracket defined below coincides on $V_i$ with that of the symmetric Lie algebra up to a positive constant \cite[proof of Lemma~4.1]{EH}.

On $L=\widehat\fg\oplus V$ define a bracket, for $\xi,\eta\in\widehat\fg$ and $u,v\in V$, by
\begin{align*}
[\xi,\eta]&=\xi\eta-\eta\xi,\\
[\xi,u]&=\xi u=-[u,\xi],\\
\tau\big([u,v],\xi\big)&=\ip{\xi u}{v}.
\end{align*}
Since $\tau$ is positive definite, the last identity determines a unique $[u,v]\in\widehat\fg$. The bracket is skew, $\sigma=\mathrm{id}_{\widehat\fg}\oplus(-\mathrm{id}_V)$ preserves it, and $[u,v,w]=[[u,v],w]$ is a trilinear product on $V$.

\begin{theorem}\label{thm:canonical-triple}
Let $(V,G,F)$ be an Eaton triple with $G$ connected, and let $V^G$ be the subspace of $G$-fixed vectors. Then:
\begin{enumerate}
\item $L=\widehat\fg\oplus V$ is a Lie algebra, and $(L,\sigma)$ is an orthogonal symmetric Lie algebra of compact type, up to the abelian factor $V^G$. Its isotropy representation is the action of $\widehat G$ on $V$. Consequently $[u,v,w]=[[u,v],w]$ is a Lie triple system on $V$.
\item $W$ is a maximal abelian subspace of $V$: $[W,W]=0$, and $[v,W]=0$ implies $v\in W$.
\item $H$ is the restricted Weyl group of $(L,\sigma)$ with respect to $W$, and $F$ is a closed Weyl chamber.
\item $L$, its bracket and the triple product depend only on the inner product of $V$ and on the partition of $V$ into $G$-orbits. More precisely, if $(V',G',F')$ is another Eaton triple with $G'$ connected, with associated $\widehat\fg'$ and $L'$, and $\varphi:V\to V'$ is a linear isometry that maps the $G$-orbits onto the $G'$-orbits, then $\xi+u\mapsto\varphi\xi\varphi^{-1}+\varphi u$ is an isomorphism of $L$ onto $L'$ which commutes with the involutions.
\end{enumerate}
Changing the sign of $[u,v]$ gives the noncompact dual $L^*$, which realizes $(V,G,F)$, up to orbit equivalence, as the isotropy representation of a symmetric space of noncompact type times the Euclidean space $V^G$, with $W$ a maximal abelian subspace of $\fp\simeq V$ and $F$ a closed restricted Weyl chamber.
\end{theorem}
\begin{proof}
(1) We first verify the Jacobi identity. By Lemma~\ref{lem:hull} there is an orthogonal decomposition
\[
V=V^G\oplus V_1\oplus\cdots\oplus V_k,
\qquad
\widehat\fg=\bigoplus_{i=1}^k\widehat\fg_i,
\]
where $\widehat\fg_i$ acts on $V_i$ as an irreducible $s$-representation and trivially on the other summands. The trace metric is orthogonal with respect to this decomposition, and its restriction to $\widehat\fg_i$ is the trace metric of the representation on $V_i$.

The bracket respects the decomposition as well: if $u\in V_i$ and $v\in V_j$ with $i\ne j$, then $\ip{\xi u}{v}=0$ for every $\xi\in\widehat\fg$, so $[u,v]=0$ by the defining identity, and the same argument gives $[V^G,V]=0$. Hence $L$ is the direct sum, as a vector space with bracket, of $V^G$ and the spaces $\widehat\fg_i\oplus V_i$.

For each $i$, Lemma~\ref{lem:hull} identifies $\widehat G_i$ with the connected isotropy group of an irreducible symmetric space, acting on $V_i$ as an $s$-representation. On $\widehat\fg_i$ the form $\tau$ is precisely the trace metric of this representation, namely the negative of the Killing form of its symmetric Lie algebra restricted to the isotropy algebra; by isotropy irreducibility, the given inner product on $V_i$ is a positive multiple of the corresponding symmetric-space metric. Consequently the bracket prescribed above agrees with the symmetric-pair bracket on $[V_i,V_i]$ up to a positive scalar. Such a rescaling preserves the Jacobi identity, exactly as in the trace-metric criterion of Eschenburg and Heintze \cite[Sections~3--4, especially Lemma~4.1]{EH}. Hence the Jacobi identity holds on every irreducible summand. The cross-summand brackets and the brackets with $V^G$ vanish, so the Jacobi identity also holds on their direct sum, and $L$ is a Lie algebra.

Consider next the inner product on $L$ which restricts to $\tau$ on $\widehat\fg$ and to the given inner product on $V$, with $\widehat\fg\perp V$. It is $\operatorname{ad}$-invariant: on $\widehat\fg$ by the invariance of the trace metric, and for $\xi\in\widehat\fg$ and $u,v\in V$ because
\[
\tau([u,v],\xi)=\ip{\xi u}{v}=-\ip{u}{\xi v},
\]
the last equality because $\xi$ is skew-symmetric. Hence $L$ is a compact Lie algebra. Since
\[
[\widehat\fg,\widehat\fg]\subset\widehat\fg,\qquad
[\widehat\fg,V]\subset V,\qquad [V,V]\subset\widehat\fg,
\]
$\sigma$ is an isometric involutive automorphism. Moreover $\widehat\fg$ is compactly embedded in $L$: for $\xi\in\widehat\fg$ one has $e^{\operatorname{ad}\xi}=\Ad(e^{\xi})\oplus e^{\xi}$, so the group generated by $e^{\operatorname{ad}\widehat\fg}$ is the image of the compact group $\widehat G$. The isotropy representation of $(L,\sigma)$ is the action of $\widehat\fg$ on $V$, which integrates to $\widehat G$.

We claim that the center of $L$ is $V^G$. If $\xi+v$ is central, then $\xi u=0$ and $[v,u]=0$ for every $u\in V$. The first condition gives $\xi=0$, because $\widehat\fg\subset\mathfrak{so}(V)$ acts faithfully; the second gives $\widehat\fg v=0$ by the defining identity, hence $v\in V^G$ since $G$ is connected and has the same orbits as $\widehat G$. Conversely, if $v\in V^G$ then $\widehat\fg v=\fg v=0$, so $[\xi,v]=0$ and $\tau([v,u],\xi)=\ip{\xi v}{u}=0$ for all $\xi\in\widehat\fg$ and $u\in V$, and $v$ is central. Thus $\widehat\fg$ meets the center trivially, and in $L=V^G\oplus\bigoplus_i(\widehat\fg_i\oplus V_i)$ the second summand is a $\sigma$-invariant ideal, compact with trivial center, hence semisimple. Therefore $(L,\sigma)$ is an orthogonal symmetric Lie algebra of compact type up to the abelian factor $V^G$, and the standard construction from a symmetric Lie algebra gives the Lie triple product $[u,v,w]=[[u,v],w]$ on $V$.

(2) If $u,v\in W$ and $\xi\in\widehat\fg$, then
\[
\tau([u,v],\xi)=\ip{\xi u}{v}=0,
\]
because $\widehat\fg u=\fg u$ ($G$ and $\widehat G$ have the same orbits) and $\fg u\perp W$ by Lemma~\ref{lem:orth}. As $\tau$ is nondegenerate, $[W,W]=0$.

For maximality, suppose that $v\in V$ satisfies $[v,w]=0$ for every $w\in W$. Then for $\xi\in\widehat\fg$ and $w\in W$
\[
0=\tau([v,w],\xi)=\ip{\xi v}{w}=-\ip{v}{\xi w},
\]
so $v$ is orthogonal to $\widehat\fg w=\fg w$. Regular points of $W$ exist by Theorem~\ref{thm:gichev-prop}, and for such $w$ the normal space is $\nu_w=W$; hence $v\in(\fg w)^\perp=\nu_w=W$.

(3) By (2), $W$ is a maximal abelian subspace. The adjoint action of $\widehat\fg$ on the summand $V$ of $L$ is its given action by skew-symmetric operators, so $e^{\operatorname{ad}\xi}$ restricts on $V$ to $e^\xi$ for $\xi\in\widehat\fg$, and these transformations generate $\widehat G$, since $\widehat G$ is connected with Lie algebra $\widehat\fg$. Therefore the restricted Weyl group is
\[
N_{\widehat G}(W)/Z_{\widehat G}(W)=H,
\]
where the equality is Lemma~\ref{lem:hull}. By Theorem~\ref{thm:main}, Corollary~\ref{cor:niezgoda} and Theorem~\ref{thm:chamber}, $F$ is the closure of a chamber of this group, that is, a closed restricted Weyl chamber.

(4) Let $\varphi$ be as in the statement. Since it maps each $G$-orbit onto a $G'$-orbit, it maps tangent spaces onto tangent spaces:
\[
\varphi(\fg x)=\fg'(\varphi x)\qquad(x\in V).
\]
By the definition of the orbit hull, this gives
\[
\varphi\widehat\fg\varphi^{-1}=\widehat\fg'.
\]
Conjugation by $\varphi$ preserves the commutator of operators and both traces in the definition of $\tau$, so
\[
\tau'(\varphi\xi\varphi^{-1},\varphi\eta\varphi^{-1})=\tau(\xi,\eta).
\]
For $u,v\in V$ and $\xi\in\widehat\fg$,
\begin{align*}
\tau'\big([\varphi u,\varphi v],\varphi\xi\varphi^{-1}\big)
&=\ip{\varphi\xi u}{\varphi v}=\ip{\xi u}{v}\\
&=\tau([u,v],\xi)
=\tau'\big(\varphi[u,v]\varphi^{-1},\varphi\xi\varphi^{-1}\big),
\end{align*}
and since $\tau'$ is nondegenerate, $[\varphi u,\varphi v]=\varphi[u,v]\varphi^{-1}$. Brackets involving elements of $\widehat\fg$ are preserved by definition of conjugation. Hence $\xi+u\mapsto\varphi\xi\varphi^{-1}+\varphi u$ is a Lie algebra isomorphism commuting with the involutions; in particular it preserves the Lie triple products.

Finally, the noncompact dual is obtained by replacing $[u,v]$ with $-[u,v]$ for $u,v\in V$ and keeping the other brackets; equivalently, by passing from $\widehat\fg\oplus V$ to $\widehat\fg\oplus iV$ inside the complexification of $L$. The isotropy action on $V$ is unchanged, as are $W$, its Weyl group and its chambers, and $V^G$ becomes the Euclidean factor.
\end{proof}

\begin{remark}[Computing the triple]
The orbit hull $\widehat\fg$ is the intersection of the subspaces $\{A\in\mathfrak{so}(V):Ax\in\fg x\}$, $x\in V$, which is attained by finitely many of them, so the linear conditions $Ax\in\fg x$ need only be imposed at finitely many suitable points. Once $\widehat\fg$ is known, $\tau$ can be computed, and $[u,v]$ is the $\tau$-dual of the functional $\xi\mapsto\ip{\xi u}{v}$: if $P$ is the Hilbert--Schmidt orthogonal projection of $\mathrm{End}(V)$ onto $\widehat\fg$ and $S$ is the positive operator on $\widehat\fg$ defined by $\tau(\xi,\eta)=\operatorname{tr}\big((S\xi)^{\top}\eta\big)$, then $[u,v]=S^{-1}P(v\otimes u)$, where $(v\otimes u)(w)=\ip{w}{u}v$. For instance, take the standard action of $G=\mathrm{SO}(2)$ on $V=\R^2$, and put $J(u_1,u_2)=(-u_2,u_1)$. Then $\widehat\fg=\R J$, $\tau(J,J)=-\operatorname{tr}(J^2)=2$, and the defining identity gives
\[
[u,v]=\tfrac12(u_1v_2-u_2v_1)J,\qquad
[u,v,w]=\tfrac12(u_1v_2-u_2v_1)Jw.
\]
The factor $1/2$ is fixed by the trace metric.

For example, the orbit hull for $G_2$ acting on $\R^7$ is $\mathfrak{so}(7)$; for $\mathrm{SO}(2)\times G_2$ acting on $\R^2\otimes\R^7$ it is $\mathfrak{so}(2)\oplus\mathfrak{so}(7)$, the isotropy algebra of the Grassmannian of $2$-planes in $\R^9$; for $\mathrm{SU}(p)\times \mathrm{SU}(q)$ acting on $\mathbb C^p\otimes\mathbb C^q$, with $p,q\ge2$ and $p\ne q$, it is the image of $\mathfrak u(p)\oplus\mathfrak u(q)$. In all these cases the bracket built from $\fg$ itself fails the Jacobi identity, since by \cite[Lemma~4.1]{EH} this requires $G$ to act as an $s$-representation, and none of these groups does: $G_2$ acts on $\R^7$ with cohomogeneity one, and the rank one symmetric spaces of dimension $7$ have constant curvature and isotropy group $\mathrm{SO}(7)$ \cite[Ch.~X]{Helgason}; $\mathrm{SO}(2)\times G_2$ and $\mathrm{SU}(p)\times\mathrm{SU}(q)$ act with cohomogeneity at least two and are properly contained in connected groups with the same orbits \cite{EH2}, while an irreducible $s$-representation of cohomogeneity at least two is maximal among the connected groups with its orbits \cite[Proposition~5.1]{EH}.
\end{remark}

\begin{remark}[Relaxing the axioms]\label{rem:relax}
Neither axiom alone is restrictive: $F=V$ satisfies (A1) for every representation, and any ray satisfies (A2). The construction above uses only (A1) and the infinitesimal form of (A2), $\ip{x}{\xi y}=0$ for $x,y\in F$ and $\xi\in\fg$ (Lemma~\ref{lem:orth}), which make $W$ a section; parts (1), (2) and (4) of Theorem~\ref{thm:canonical-triple} use nothing else. These two conditions do not force $F$ to be a chamber: $F=W$, or the union of two adjacent closed chambers for $\mathrm{SU}(3)$ acting on the traceless Hermitian $3\times3$ matrices, satisfies both but not (A2). The difference is that between critical points and maxima: for $x,y\in W$ regular, the critical points of $\ip{x}{\cdot}$ on $Gy$ are the $|H|$ points of $Hy$, one in each chamber, and (A2) asks that the one in the chamber of $x$ be the maximum. By Theorems~\ref{thm:chamber} and \ref{thm:char}, a closed convex cone in $W$ meeting every orbit satisfies (A2) if and only if it is a closed Weyl chamber. Thus the infinitesimal conditions determine the Lie triple system, while (A2) selects a chamber.
\end{remark}

\subsection{Uniqueness and non-canonical structures}\label{sec:noncanonical}

\begin{proposition}[Other Lie triple structures]\label{prop:other-triples}
Let $K\subset \mathrm{O}(V)$ be a compact connected subgroup with the same orbits as $G$ which acts on $(V^G)^\perp$ as an $s$-representation (the isotropy representation of a symmetric space of compact or, equivalently, of noncompact type), and let $L_K=\mathrm{Lie}(K)\oplus V$, with the bracket defined as above using the trace metric of $(\mathrm{Lie}(K),V)$. Then $L_K$ is an orthogonal symmetric Lie algebra, $W$ is a maximal abelian subspace of $V$ for it, its restricted Weyl group is $H$, and $F$ is a closed Weyl chamber. Moreover $K\subset\widehat G$, and $L_K=L$, as Lie algebra structures on $\mathfrak{so}(V)\oplus V$, if and only if $K=\widehat G$.
\end{proposition}
\begin{proof}
Groups with the same orbits have the same fixed vectors and invariant subspaces. Thus $V^K=V^G$, and the $G$-irreducible summands $V_i$ of $(V^G)^\perp$ are $K$-irreducible, hence the tangent spaces of the irreducible factors of the symmetric space, and $\mathrm{Lie}(K)=\bigoplus_i\mathfrak k_i$, where the isotropy algebra $\mathfrak k_i$ of the $i$-th factor acts on $V_i$ as an irreducible $s$-representation and trivially on the other summands. The proof of Theorem~\ref{thm:canonical-triple}(1), with \cite[Lemma 4.1]{EH} applied to each pair $(\mathfrak k_i,V_i)$, shows that $L_K$ is a Lie algebra. Since the orbits agree, $\mathrm{Lie}(K)u=\fg u$ for every $u$, and $(V,K,F)$ is an Eaton triple, exactly as for $\widehat G$ in the proof of Lemma~\ref{lem:hull}. Theorem~\ref{thm:main}, applied to both triples, gives $N_K(W)x=Kx\cap W=Gx\cap W=N_G(W)x$ for $x\in W$, and the finite-group argument in the proof of Theorem~\ref{thm:connected}, which uses only that both groups are finite with the same orbits on $W$, gives $N_K(W)/Z_K(W)=H$ as subgroups of $\mathrm{O}(W)$. The proofs of Theorem~\ref{thm:canonical-triple}(2) and (3) now apply with $K$ in place of $\widehat G$ and this equality in place of Lemma~\ref{lem:hull}. Finally, $K$ preserves the $G$-orbits, so $K\subset\widehat G$ by Lemma~\ref{lem:hull}; if $K\neq\widehat G$ then $\mathrm{Lie}(K)\neq\widehat\fg$, since both groups are connected, and $L_K\neq L$.
\end{proof}

By Dadok's theorem such groups $K$ always exist, and each realizes $(V,G,F)$, up to orbit equivalence, as the isotropy representation of a symmetric space. They need not be unique.

\begin{remark}[Non-canonical structures]\label{rem:noncanonical-lie-triple}
(a) \emph{Spheres.} Take $V=\R^{2m}$, $m\ge2$, and $F=\R_{\ge0}e_1$. Both $\mathrm{SO}(2m)$ and $\mathrm{U}(m)$ have the spheres as orbits, and both are isotropy representations of symmetric spaces, namely $\mathrm{SO}_0(2m,1)/\mathrm{SO}(2m)$ and $\mathrm{SU}(m,1)/\mathrm{S}(\mathrm{U}(m)\times \mathrm{U}(1))$. The corresponding Lie triple systems are not isomorphic. Indeed, an isomorphism of Lie triple systems $\varphi:V\to V'$ satisfies $\varphi[u,v,\cdot]\varphi^{-1}=[\varphi u,\varphi v,\cdot]$, hence extends to an isomorphism of the Lie algebras $[V,V]\oplus V$, where $[V,V]$ is the span of the operators $[u,v,\cdot]$ with the commutator bracket. Here these are $\mathfrak{so}(2m+1)$ and $\mathfrak{su}(m+1)$, dually the transvection algebras $\mathfrak{so}(2m,1)$ and $\mathfrak{su}(m,1)$ of the two hyperbolic spaces, of dimensions $m(2m+1)$ and $m(m+2)$, which differ for $m\ge2$. Geometrically, real hyperbolic space has constant curvature and complex hyperbolic space does not. On $\R^{4n}$, $n\ge2$, there is a third structure, from $\mathrm{Sp}(n)\cdot \mathrm{Sp}(1)$ and quaternionic hyperbolic space, and on $\R^{16}$ a fourth, from $\mathrm{Spin}(9)$ and the Cayley hyperbolic plane; the same dimension count, with $\mathfrak{sp}(n+1)$ and $\mathfrak f_4$, distinguishes them. All these groups act transitively on the spheres \cite[Introduction]{EH2}, so all give Eaton triples with the same cone $F$ and the same group $H=\{\pm1\}$, of type $A_1$.

(b) \emph{Why they are not canonical.} The additional complex, quaternionic and octonionic structures are not determined by the orbit partition of $V$, so different isotropy representations may have the same orbits and different Lie triple systems. The largest connected group with these orbits is $\widehat G=\mathrm{SO}(2m)$, respectively $\mathrm{SO}(4n)$ and $\mathrm{SO}(16)$, and gives the triple system of the sphere, dually of real hyperbolic space.

(c) \emph{Uniqueness within the trace-metric construction.} If $V$ is $G$-irreducible of cohomogeneity at least two, an irreducible $s$-representation with the orbits of $G$ is maximal among connected groups with these orbits \cite[Proposition~5.1]{EH}. Hence every $K$ as in Proposition~\ref{prop:other-triples} equals $\widehat G$ and $L_K=L$: the trace-metric construction gives the only Lie triple structure arising from connected orbit-equivalent isotropy groups with their trace metrics, and when $G$ itself acts as an $s$-representation, the construction from $G$ gives this structure. In the reducible case the same holds summand by summand, so different structures can arise only on summands of cohomogeneity one, as in (a). Rescaling $[V_i,V_i]$ by a positive constant on an irreducible summand still gives a compact symmetric Lie algebra with the same isotropy action; it is the trace metric that fixes the normalization.
\end{remark}

For disconnected $G$ no such structure exists in general: in Example~\ref{ex:I25} the group $H=I_2(5)$ is not crystallographic, while the restricted Weyl group of a symmetric Lie algebra always is.

\subsection{Eaton triples and Lie groups with a bi-invariant distance}\label{sec:biinvariant}

In the adjoint case Eaton triples meet the geometry of Lie groups with a bi-invariant metric, studied by Milnor \cite{Milnor} in the Riemannian case and in \cite{LR} for Finsler metrics.

\begin{remark}[The group case]\label{rem:group-case}
For a compact connected Lie group $K$ with Lie algebra $\mathfrak k$ and the Eaton triple $(\mathfrak k,\Ad K,F)$, $F$ a closed Weyl chamber, the adjoint representation of the semisimple part is already maximal: on each simple ideal it is an irreducible $s$-representation, either of rank at least two, hence maximal by \cite[Proposition~5.1]{EH}, or the full group $\Ad \mathrm{SU}(2)=\mathrm{SO}(3)$. Hence $L$ is the direct sum of the compact symmetric Lie algebras $\mathfrak k_j\oplus\mathfrak k_j$ of the simple ideals and of the center of $\mathfrak k$, and on each simple ideal the triple product is a positive multiple of $[[u,v],w]$. Thus in the group case the construction returns the triple product of $\mathfrak k$, up to scaling on simple ideals, while in general it returns a Lie triple product on $V$ which need not come from a bracket on $V$.
\end{remark}

In this subsection $G$ is a Lie group with Lie algebra $\fg$, neither compact nor connected a priori, $\|\cdot\|$ is an $\Ad$-invariant norm on $\fg$, and we consider Eaton triples of the form $(\fg,\overline{\Ad G},F)$. Being a closed bounded subgroup of the isometry group of the norm, $\overline{\Ad G}$ is compact, and averaging gives an $\Ad G$-invariant inner product on $\fg$; hence $\fg=\mathfrak z\oplus[\fg,\fg]$, with $\mathfrak z$ the center and $[\fg,\fg]$ compact semisimple, as in \cite{Milnor}, see also \cite[Corollary 4.25 and its proof]{Knapp}.

Assume now that $G$ is connected, let $\mathfrak t$ be a maximal abelian subalgebra of $[\fg,\fg]$ and $C\subset\mathfrak t$ a closed Weyl chamber. Then
\[
\big(\fg,\ \overline{\Ad G},\ F\big),\qquad F=\mathfrak z\oplus C,
\]
is an Eaton triple, with reduced triple $(\mathfrak z\oplus\mathfrak t,\,\mathcal W,\,F)$, where $\mathcal W$ is the Weyl group of $[\fg,\fg]$ acting trivially on $\mathfrak z$; the fixed vectors are exactly $\mathfrak z$, so $\mathcal W$ is not essential relative to $\mathfrak z\oplus\mathfrak t$ in the terminology of \cite{DolberryTam}. Since $\Ad G$ fixes $\mathfrak z$ and $\mathfrak z\perp[\fg,\fg]$, it suffices to check (A2) for $x,y\in C$, and this follows from two results of Kostant: for $y\in C$, the orthogonal projection of $\Ad(g)y$ onto $\mathfrak t$ lies in the convex hull of $\mathcal W y$ \cite[Theorem 8.2]{Kostant}, applied to the complexification of $[\fg,\fg]$ regarded as a real Lie algebra with Cartan decomposition $[\fg,\fg]\oplus i[\fg,\fg]$, and a linear form given by an element of $C$ attains its maximum on $\mathcal W y$ at $y$ \cite[Lemma 3.2]{Kostant}. We remark that (A2) is a consequence of the convexity theorem, not a restatement of it.

If $G$ is disconnected but $G/G^0$ is finite, then $\Ad G$ is closed, $\Ad G^0$ fixes $\mathfrak z$, and the action of $\Ad G$ is still polar with section $\mathfrak z\oplus\mathfrak t$; by Theorem~\ref{thm:char} it admits an Eaton cone if and only if the group induced on $\mathfrak z\oplus\mathfrak t$ by its normalizer in $\Ad G$, which contains $\mathcal W$, is a finite reflection group, and the cone is then a closed chamber of this larger group.

If $G/G^0$ is infinite, $\mathfrak z\oplus\mathfrak t$ need not be a section, since the identity component of $\overline{\Ad G}$ may rotate the center: for $G=\R^2\rtimes\Z$, with the generator acting by a rotation of angle $2\pi\theta$, $\theta$ irrational, and the Euclidean norm, $\fg=\mathfrak z=\R^2$ and $\overline{\Ad G}=\mathrm{SO}(2)$, whose sections are the lines through the origin.

For connected $G$, restriction is a bijection between the $\Ad G$-invariant norms on $\fg$ and the $\mathcal W$-invariant norms on $\mathfrak z\oplus\mathfrak t$ \cite[Theorem 6]{TamHill}, \cite[Theorem 6.5]{Lewis}; both results assume the existence of the reduced triple (a subsystem in the sense of \cite[Assumption 4.1]{Lewis}), which Theorem~\ref{thm:main} provides, and the adjoint case is discussed in \cite[Example 3]{TamHill}. If $G$ is disconnected and $\Ad G$ admits an Eaton cone, the same holds with the group of the full reduced triple in place of $\mathcal W$; the distinction matters because the components other than $G^0$ may also act on $\mathfrak z$. Smoothness, strict convexity, the extreme and exposed points of the unit ball \cite[Theorem 6.5]{Lewis}, and norming functionals and faces \cite{DolberryTam} can also be studied through the reduced triple. These are the $\Ad$-invariant norms whose metric geometry is studied in \cite{LR}.

On the other hand, the triple depends only on $\fg$ and $\overline{\Ad G}$, so it forgets the global form of $G$: $G$ and its coverings, or $K\times\R^k$ and $K\times T^k$ for $K$ compact, give the same triple, and the abelian factor in the structure theorem of \cite{LR} survives only as the subspace $\mathfrak z$ of fixed vectors. The triple describes all invariant norms but singles out none; in the direction of $\mathfrak z$ invariance imposes no condition, so a norm need not split along $\mathfrak z\oplus[\fg,\fg]$. Finally, length, distance, short paths and flat planes, the subject of \cite{LR}, are not orbit data.

Conversely, not every Eaton triple comes from a Lie group. By Theorem~\ref{thm:canonical-triple}, those with connected group are, up to orbit equivalence, isotropy representations of symmetric spaces, and the adjoint ones correspond to symmetric spaces of group type $(K\times K)/\Delta K$, times a Euclidean factor. The spheres of Remark~\ref{rem:noncanonical-lie-triple} and the spaces of $p\times q$ matrices are in general not of this type, and Example~\ref{ex:I25} cannot arise from the isotropy representation of a symmetric space, since its reduced group is noncrystallographic. A common framework, and a setting for further study, is that of symmetric spaces $M=U/K$ with an invariant Finsler metric, which correspond to the $\Ad K$-invariant norms on the Eaton triple at the base point; Lie groups with a bi-invariant metric are the case of group type.

\par\bigskip
\noindent\textbf{Acknowledgement.}\ 
Most of this research was conducted during Gabriel Larotonda's visit to the Department of Mathematics and Statistics at the University of Nevada, Reno, in October 2026. Gabriel Larotonda thanks Tin-Yau Tam for his hospitality.

\par\bigskip
\noindent\textbf{Funding.}\ 
Gabriel Larotonda was supported by Universidad de Buenos Aires, grant UBACyT 20020220400256BA.

\par\bigskip
\noindent\textbf{AI disclosure.}\ 
In preparing this manuscript, Gabriel Larotonda used Claude (Anthropic) and Tin-Yau Tam used ChatGPT (OpenAI) for assistance with exposition, organization, and the review of mathematical arguments. The authors are responsible for the mathematical content, references, and final text.

\end{document}